\documentclass[nospthms,envcountsect]{svmult} 
\usepackage{latexsym}
\usepackage{amsfonts}
\usepackage{amsmath}
\usepackage{amssymb}
\usepackage{amscd}
\usepackage{enumerate}
    \newtheorem{theorem}                    {Theorem}       [section]
    \newtheorem{lemma}      [theorem]       {Lemma}
    
    \newtheorem{proposition}[theorem]       {Proposition}

\begin{document}
\catcode`@=11
\atdef@ I#1I#2I{\CD@check{I..I..I}{\llap{$\m@th\vcenter{\hbox
  {$\scriptstyle#1$}}$}
  \rlap{$\m@th\vcenter{\hbox{$\scriptstyle#2$}}$}&&}}
\atdef@ E#1E#2E{\ampersand@
  \ifCD@ \global\bigaw@\minCDarrowwidth \else \global\bigaw@\minaw@ \fi
  \setboxz@h{$\m@th\scriptstyle\;\;{#1}\;$}%
  \ifdim\wdz@>\bigaw@ \global\bigaw@\wdz@ \fi
  \@ifnotempty{#2}{\setbox@ne\hbox{$\m@th\scriptstyle\;\;{#2}\;$}%
    \ifdim\wd@ne>\bigaw@ \global\bigaw@\wd@ne \fi}%
  \ifCD@\enskip\fi
    \mathrel{\mathop{\hbox to\bigaw@{}}%
      \limits^{#1}\@ifnotempty{#2}{_{#2}}}%
  \ifCD@\enskip\fi \ampersand@}
\catcode`@=\active

\renewcommand{\labelenumi}{\alph{enumi})}
\newcommand{\chr}{\operatorname{char}}
\newcommand{\isom}{\stackrel{\sim}{\longrightarrow}}
\newcommand{\Aut}{\operatorname{Aut}}
\newcommand{\Shv}{\operatorname{Shv}}
\newcommand{\Hom}{\operatorname{Hom}}
\newcommand{\End}{\operatorname{End}}
\newcommand{\HOM}{\operatorname{{\mathcal H{\frak{om}}}}}
\newcommand{\Mod}{\operatorname{Mod}}
\newcommand{\EXT}{\operatorname{\mathcal E{\frak xt}}}
\newcommand{\Tot}{\operatorname{Tot}}
\newcommand{\Ext}{\operatorname{Ext}}
\newcommand{\Gal}{\operatorname{Gal}}
\newcommand{\cosk}{\operatorname{cosk}}
\newcommand{\Pic}{\operatorname{Pic}}
\newcommand{\Spec}{\operatorname{Spec}}
\newcommand{\trdeg}{\operatorname{trdeg}}
\newcommand{\holim}{\operatornamewithlimits{holim}}
\newcommand{\im}{\operatorname{im}}
\newcommand{\hyp}{{\rm hyp}}
\newcommand{\coim}{\operatorname{coim}}
\newcommand{\coker}{\operatorname{coker}}
\newcommand{\gr}{\operatorname{gr}}
\newcommand{\id}{\operatorname{id}}
\newcommand{\Br}{\operatorname{Br}}
\newcommand{\cd}{\operatorname{cd}}
\newcommand{\CH}{CH}
\newcommand{\Alb}{\operatorname{Alb}}
\renewcommand{\lim}{\operatornamewithlimits{lim}}
\newcommand{\colim}{\operatornamewithlimits{colim}}
\newcommand{\rk}{\operatorname{rank}}
\newcommand{\codim}{\operatorname{codim}}
\newcommand{\NS}{\operatorname{NS}}
\newcommand{\cone}{{\rm cone}}
\newcommand{\rank}{\operatorname{rank}}
\newcommand{\ord}{{\rm ord}}
\newcommand{\f}{{\cal F}}
\newcommand{\g}{{\cal G}}
\newcommand{\p}{{\cal P}}
\newcommand{\N}{{\mathbb N}}
\newcommand{\A}{{\mathbb A}}
\newcommand{\Z}{{{\mathbb Z}}}
\newcommand{\Q}{{{\mathbb Q}}}
\newcommand{\R}{{{\mathbb R}}}
\newcommand{\B}{{\mathbb Z}^c}
\renewcommand{\H}{{{\mathbb H}}}
\renewcommand{\P}{{{\mathbb P}}}
\newcommand{\F}{{{\mathbb F}}}
\newcommand{\m}{{\mathfrak m}}
\newcommand{\Sm}{{\text{\rm Sm}}}
\newcommand{\Sch}{{\text{\rm Sch}}}
\newcommand{\et}{{\text{\rm et}}}
\newcommand{\eh}{{\text{\rm eh}}}
\newcommand{\Zar}{{\text{\rm Zar}}}
\newcommand{\Nis}{{\text{\rm Nis}}}
\newcommand{\tr}{\operatorname{tr}}
\newcommand{\tor}{{\text{\rm tor}}}
\newcommand{\PreShv}{\text{\rm PreShv}}
\newcommand{\Div}{\operatorname{Div}}
\newcommand{\Ab}{{\text{\rm Ab}}}
\renewcommand{\div}{\operatorname{div}}
\newcommand{\corank}{\operatorname{corank}}
\renewcommand{\O}{{\cal O}}
\newcommand{\C}{{\cal C}}
\renewcommand{\p}{{\mathfrak p}}
\newcommand{\proof}{\noindent{\it Proof. }}
\newcommand{\proofend}{\hfill $\Box$ \\}
\newcommand{\rem}{\noindent {\it Remark. }}
\newcommand{\example}{\noindent {\bf Example. }}
\newcommand{\ar}{{\text{\rm ar}}}
\newcommand{\del}{{\delta}}

\title*{Parshin's conjecture and motivic cohomology with compact support}
\author{Thomas Geisser\thanks{Supported in part by NSF grant No.0556263}}
\institute{University of Southern California}

\maketitle

\begin{abstract} We discuss Parshin's conjecture on rational $K$-theory over
finite fields and its implications for motivic cohomology
with compact support.
\end{abstract}

\section{Introduction}
Parshin's conjecture states that higher algebraic $K$-groups of smooth
projective schemes over finite fields are torsion.
In \cite{ichparshin}, we studied the properties that Parshin's 
conjecture would imply for rational higher Chow groups. We compared
higher Chow groups to weight homology $H_i^W(X,\Q(n))$, 
defined by Jannsen \cite{jannsen}
based on the work of Gillet-Soule \cite{gilletsoule}, and obtained
a diagram
\begin{equation}\label{maina}
\begin{CD}
 H_i^c(X,\Q(n)) @>\pi>>  H_i^W(X,\Q(n))\\
@V\alpha VV @A\gamma AA  \\
\tilde H_i^c(X,\Q(n))@>\beta >>\tilde  H_i^W(X,\Q(n)).
\end{CD}\end{equation}
The terms with the tilde are the cohomology of the first
non-vanishing $E^1$-line of the niveau spectral sequence.
Parshin's conjecture in weight $n$ is equivalent to $\pi$
being an isomorphism for all $X$ and $i$. We showed that
$\pi$ is an isomorphism if and only if $\alpha$, $\beta$
and $\gamma $ are isomorphisms, and gave criteria for 
this to happen.

In this article, we take the cohomological point of view and
examine the properties that Parshin's conjecture implies
for motivic cohomology with compact support. Surprisingly,
the properties obtained are not dual to the properties
for higher Chow groups, but have a different flavor. 
The method to study motivic cohomology with compact support
is to use the coniveau filtration. To avoid the problems
arising from the covariance of motivic cohomology with 
compact support for open embeddings (for example, one gets very large by
taking inverse limits, and has to deal with derived inverse
limits), we consider the dual groups $H^i_c(X,\Q(n))^*$. We obtain
a niveau spectral sequence, and compare it with the  
spectral sequence for the dual of weight cohomology 
$H^i_W(X,\Q(n))^*$ as in \cite{ichparshin} to obtain a diagram 
\begin{equation}\label{mainb}
\begin{CD}
\tilde H^i_W(X,\Q(n))^* @= \tilde H^i_c(X,\Q(n))^*\\
@V\gamma^* VV @V\alpha^*VV  \\
H^i_W(X,\Q(n))^*@>\pi^*   >>  H^i_c(X,\Q(n))^*.
\end{CD}\end{equation}
Again, the upper terms are given by the first non-vanishing row
of $E^1$-terms in the niveau spectral sequence. 
The map $\pi^*$ is an isomorphism for all $X$ if and only if
Parshin's conjecture holds.  
In contrast to the homological situation,
$\alpha^*$ being an isomorphism is stronger than Parshin's conjecture.
We go on to examine the relationship between diagrams 
\eqref{maina} and \eqref{mainb}. Not surprisingly, this is related to 
Beilinson's conjecture that rational and numerical equivalence agree up
to torsion over finite fields.
Finally we relate bounds for all four rational motivic theories
to Parshin's conjecture.   

Since the purpose of this work is to understand interrelations between
certain conjectures, we assume the existence of resolution of 
singularities. Its use in the results of Friendlander and Voevoesky 
\cite{friedvoe} maybe be dispensable with more work because
we work with rational coefficients, 
but occasionally we need a smooth and projective model for every 
function field to do an induction process. 

Throughout this paper, the cateogory of schemes over $k$, written
$Sch/k$ denotes the cateogory of separated schemes of finite
type over $k$, and $Sm/k$ the category of smooth schemes over $k$.

\medskip
{\it Acknowledgements:}  This paper was inspired by the work of and
discussions with U.Jannsen and S.Saito.

\section{Motivic cohomology with compact support}

For a scheme $X$ over a field $k$, motivic cohomology with
compact support is defined as
$$ H^i_c(X,\Z(n))=\Hom_{DM^-}(M^c(X),\Z(n)[i]).$$
A concrete description is given as follows \cite[\S 3]{friedvoe}:
Let $\rho:(Sch/k)_{cdh}\to (Sm/k)_{Nis}$ be the map from the 
large cdh-site of $k$ to the smooth site with the Nisnevich topology.
Let $\Z(n)$ be the motivic complex on $(Sm/k)_{Nis}$, and consider
an injective resolution $\rho^* \Z(n)\to I^\cdot$
on $(Sch/k)_{cdh}$ (we need resolution of singularites to 
ensure that $\rho^*$ is exact). Let $\Z^c(X)$
be the cdh-sheafification of the presheaf which associates to 
$U$ the free abelian group
generated by those subschemes $Z\subseteq X\times U$ whose 
projection to $U$ induces an open embedding $Z\subseteq U$. 
Then $H^i_c(X,\Z(n))=\Hom_{D(Shv_{cdh})}(\Z^c(X),I^\cdot[i])$.
This satisfies the following properties:
\begin{enumerate}
\item Contravariance for proper maps.
\item Covariance for flat quasi-finite maps.
\item For a closed subscheme $Z$ of $X$ with open complement $U$,
there is a localization sequence
\begin{equation}\label{locali}
\cdots  \to H^i_c(U,\Z(n))\to  H^i_c(X,\Z(n))\to  H^i_c(Z,\Z(n)) \to\cdots.
\end{equation}
\end{enumerate}
If $X$ is proper, then since $\Z^c(X)=\Z(X)$, motivic cohomology with 
compact support agrees with motivic cohomology 
$H^i(X,\Z(n)):= H^i_{cdh}(X,\Z(n))$. Moreover, 
under resolution of singularities, we get for smooth $X$ 
of dimension $d$ isomorphisms \cite{susvoe00, voevodskysmooth}
\begin{equation}\label{smooth}
H^i_{cdh}(X,\Z(n)) \cong H^i_{Nis}(X,\Z(n))\cong CH^n(X,2n-i).
\end{equation}

\begin{proposition}\label{vanbound}
a) We have $H^i_c(X,\Z(n))=0$ for $i>n+\dim X$. 

b) If $k$ is finite, resolution of singularities exists,
and if $n>\dim X$, then $H^i_c(X,\Q(n))=0$ for $i\geq n+\dim X$.

c) If $k$ is finite and $X$ is smooth of dimension $d$, then
$H^{n+d}(X,\Q(n))=0$ unless $n=d$.
\end{proposition}

\proof
a) Using the localization sequence and induction on the dimension, 
the statement is easily reduced to the case where $X$ is proper.
Then we use that the complex $\Z(n)$ is concentrated in degrees at most 
$n$, and $X$ has $cdh$-cohomological dimension $d$.

b) This was proved in \cite[Prop.6.3]{ichsuslin}. The idea is to 
use induction on the dimension to reduce to $X$ smooth and proper,
and then use c).

c) If $n<d$ then this follows by comparing to higher Chow groups.
If $n>d$, consider the spectral sequence
\begin{equation}
E_1^{s,t}=\bigoplus_{x\in X^{(s)}}H^{t-s}(k(x),\Z(n-s))
\Rightarrow H^{s+t}(X,\Z(n)).
\end{equation}
In order for the $E_1^{s,t}$-terms not to vanish, we need $t\leq n$
and $s\leq d$, hence to have $s+t=n+d$ we need $s=d$ and $t=n$.
But $E_1^{d,n}$ is a sum of $H^{n-d}(k(x),\Z(n-d))$ for finite fields $k(x)$,
and higher Milnor K-theory of finite fields is torsion.
\proofend

\subsection{The niveau spectral sequence}
In order not to deal with derived inverse limits and to get smaller
groups, we work with the dual of motivic cohomology with compact
support
$$ H^i_c(X,\Q(n))^* := \Hom(H^i_c(X,\Z(n)),\Q).$$
These groups are covariant for proper maps and contravariant
for quasi-finite flat maps. 
Let $Z_s$ be set of closed subschemes of dimension at most $s$
and let $Z_s/Z_{s-1}$ be the set of ordered pairs
$(Z,Z')\in Z_s\times Z_{s-1}$ such that $Z'\subseteq Z$.
Then $Z_s$ as well as $Z_s/Z_{s-1}$ are ordered by inclusion,
and we obtain a filtration $Z_0\subseteq Z_1\subseteq \cdots$. 
We use covariance for proper maps to define
$$ H^i_c(Z_s,\Q(n))^*:=\colim_{Z\in Z_s}H^i_c(Z,\Q(n))^*.$$
For a point $x\in X $ we write 
$x\in Z_s$ if $\overline{\{x\}}\in Z_s$, and using contravariance
for open embeddings define
$$ H^i_c(k(x),\Q(n))^*:=\colim_{U\cap\overline{\{x\}}\not=\emptyset}
H^i_c(U\cap\overline{\{x\}},\Q(n))^*.$$
Beware that this is typically not the dual of any group. For example, 
for the function field $k(C)$ of a smooth and proper curve $C$ 
we have
$\lim_U H^1_c(U,\Q(0)) =(\prod_{C_{(0)}} \Q) /\Q$, whereas taking duals
allows us to work with the countable "predual" group 
$H^1_c(k(C),\Q(0))^*=\colim_U H^1_c(U,\Q(0))^*
=\ker\big(  \oplus_{C_{(0)}} \Q\to \Q\big) $.
From the localization sequence we obtain
$$H^i_c(Z_s/Z_{s-1},\Q(n))^*:=\colim_{(Z,Z')\in Z_s/Z_{s-1}}H^i_c(Z-Z',\Q(n)^*
=\bigoplus_{x\in Z_d}H^i_c(k(x),\Q(n))^*.$$

The usual yoga with exact couples gives

\begin{proposition}
There is a homological spectral sequence
\begin{equation}\label{sseq}
E^1_{s,t}=
\bigoplus_{x\in X_{(s)}}H^{s+t}_c(k(x),\Q(n))^*
\Rightarrow H^{s+t}_c(X,\Q(n))^*.
\end{equation}
\end{proposition}

The $d^1$-differential is given by
$$ H^{i+1}_c(Z_{s+1}/Z_s,\Q(n))^*\to H^i_c(Z_s,\Q(n))^*
\to H^i_c(Z_s/Z_{s-1},\Q(n))^*.$$
By Proposition \ref{vanbound}a), we obtain 
$H^i_c(k,\Q(n))^*=0$ for $i>n+s$, i.e. $E^1_{s,t}$ vanishes
for $t>n$, so that the spectral sequence \eqref{sseq} is concentrated below
and on the line $t=n$. 
On the line $t=n$, the terms $E^1_{s,n}$
vanish for $s<n$ by Proposition \ref{vanbound}b). 
We define $\tilde H^j_c(X,\Q(n))^*$ to be the cohomology of the line 
$E^1_{*,n}$ 
\begin{equation}\label{comp}
\bigoplus_{x\in X_{(n)}}H^{2n}_c(k(x),\Q(n))^* \leftarrow \cdots 
\leftarrow\bigoplus_{x\in X_{(d)}}H^{n+d}_c(k(x),\Q(n))^*,
\end{equation}
where we put the term indexed by $X_{(i)}$ in degree $n+i$.
It is easy to check that we obtain canonical maps
\begin{equation}\label{canmap}
\tilde H^i_c(X,\Q(n))^*\stackrel{\alpha^*}{\to}  H^i_c(X,\Q(n))^* 
\end{equation}

\section{Parshin's conjecture}
Parshin's conjecture states that for all smooth and projective
$X$ over $\F_q$, the groups $K_i(X)_\Q$ are torsion for $i>0$.
In \cite{ichtate} we showed that it is implies by Tate's conjecture
and Beilinson's conjecture that rational and numerical equivalence
agree up to torsion.
Since $K_i(X)_\Q=\oplus_n H^{2n-i}(X,\Q(n))$, it follows that Parshin's
conjecture is equivalent to the following conjecture for all $n$.

\medskip

\noindent{\bf Conjecture $P^n$:}
{\it For all smooth and projective schemes $X$
over the finite field $\F_q$, and all $i\not=2n$, the group $H^i(X,\Z(n))$
is torsion.}

\medskip

Conjecture $P^n$ is known for $n=0,1$ and is trivial for $n<0$.
In \cite{ichparshin}, we considered the homological analog
(it was denoted $P(m)$ in loc.cit.):

\medskip 

\noindent{\bf Conjecture $P_m$:}
{\it For all smooth and projective schemes $X$
over the finite field $\F_q$, and all $i\not=2m$, the group $H_i^c(X,\Z(m))$
is torsion.}

\medskip

This conjecture is not known for any $m$. 
One can also consider the restrictions $P^n(d)$ and $P_m(d)$
of the above conjectures to varieties  of dimension at most $d$. 
By the projective bundle formula one gets $P^n(d)\Rightarrow P^{n-1}(d-1)$
and $P_m(d)\Rightarrow P_{m-1}(d-1)$, hence $P^n\Rightarrow P^{n-1}$
and $P_m\Rightarrow P_{m-1}$.

\begin{lemma}
We have $P^n(d)\Leftrightarrow P_{d-n}(d)$.
\end{lemma}

\proof
Let $X$ be smooth and projective of dimension $e\leq d$. Then
conjecture $P^{n-d+e}$ holds for $X$, hence the formula
$H^i(X,\Z(a))\cong H_{2e-i}^c(X,\Z(e-a))$ implies conjecture $P_{d-n}$
for $X$. The converse is proved the same way.
\proofend

Since conjecture $P^{-1}$ is trivially true, the following Lemma
explains why the spectral sequence for homology with compact support
in \cite{ichparshin} is concentrated in degrees $s\geq n$, 
whereas \eqref{sseq} a priori is not:

\begin{lemma}
If conjecture $P_{-1}$ holds, $H^i_c(X,\Q(n))=0$ for $n>d=\dim X$ and any $X$.
In particular, the terms $E^1_{s,t}$ vanish for $s<n$ in the spectral 
sequence \eqref{sseq}.
\end{lemma}

\proof
By induction on the dimension of $X$ and the sequence \eqref{locali} we can
assume that $X$ is smooth and proper. Then 
$H^i_c(X,\Q(n))=H_{2d-i}^c(X,\Q(d-n))$ which vanishes
by conjecture $P_{-1}$.
\proofend

\begin{lemma}\label{Ha}
The following statements are equivalent:
\begin{enumerate}
\item Conjecture $P^n$.
\item For all schemes $X$ over $\F_q$, we have $H^i_c(X,\Q(n))=0$ for $i<2n$.
\item For all finitely generated fields $k/\F_q$,
we have $H^i_c(k,\Q(n))^*=0$ for $i<2n$.
\end{enumerate}
\end{lemma}

\proof
a) $\Rightarrow$ b) follows by induction on the dimension and 
localization to recude to the smooth and proper case, 
b) $\Rightarrow$ c) follows by taking colimits, and 
c) $\Rightarrow$ a) follows with the spectral sequence \eqref{sseq}.
\proofend

It is not a priori clear if the terms $H^i_c(k(x),\Q(n))$ with 
$2n\leq i<\trdeg k(x)+n$ should vanish or not. Thus the following 
statement is possibly stronger than Parshin's conjecture
(but see Proposition \ref{compareto}):

\begin{proposition}\label{themapg}
The following statements are equivalent:
\begin{enumerate}
\item Conjecture $P^n$ holds, and for smooth and projective $X$ we have
$$ \tilde H^i_c(X,\Q(n))^* \cong \begin{cases}
CH^n(X)^* &i=2n;\\
0 &\text{else}.
\end{cases}$$
\item The groups $H^i_c(k,\Q(n))^*$ vanish for $i\not=n+\trdeg k$.
\item The map $\alpha^*$ is an isomorphism for all $X$.
\end{enumerate}
\end{proposition}

\proof
a) $\Rightarrow$ b): We proceed by induction on the transcendence
degree. Choose a smooth and projective model $X$ of $k$.
Since $H^i_c(X,\Q(n))$ is $CH^n(X)$ for $i=2n$ and
vanishies for $i\not=2n$, an inspection of
the spectral sequence \eqref{sseq} shows the vanishing.
b) $\Rightarrow$ c) is clear.

c) $\Rightarrow$ a): Conjecture $P^n$ follows because 
$\tilde H^i_c(X,\Q(n))^*$ vanishes for $i< 2n$, and the sequence 
is exact because for smooth and proper $X$, 
$H^i_c(X,\Q(n))^*$ vanishes for $i>2n$ 
and is isomorphic to $CH^n(X)$ for $i=2n$.
\proofend

The statements of this Proposition are non-trivial
even in the case $n=0$ (but they can be proven with methods similar
to \cite[Thm.5.10]{jannsen} in this case).

\section{Weight cohomology}
Let $\C$ be category of correspondences with objects smooth projective 
varieties $[X]$ over the field $k$,  
$\Hom_\C([X],[Y])=\oplus CH_{\dim X_i}(X_i\times Y)_\Q$ for $X=\coprod X_i$
the decomposition into connected components, and the usual composition
of correspondenes. In \cite{gilletsoule},
Gillet and Soul\'e defined, for every separated scheme of finite type, a
weight complex $W(X)$ in the homotopy category of bounded 
complexes in $\C$, satisfying the following properties
\cite[Thm. 2]{gilletsoule}:

\begin{enumerate}
\item $W(X)$ is represented by a bounded complex
$$ [X_0]\gets [X_1]\gets \dots \gets [X_k]$$
with $\dim X_i\leq \dim X-i$.
\item $W(-)$ is covariant functorial for proper maps.
\item $W(-)$ is contravariant functorial for open embeddings.
\item If $T\to X$ is a closed embedding with open complement $U$,
then there is a distinguished triangle
$$ W(T)\stackrel{i_*}{\longrightarrow} W(X) \stackrel{j^*}{\longrightarrow}
W(U).$$
\end{enumerate}

Our notation differs from loc.cit. in variance. In loc.cit.,
resolution of singularities is used to obtain an integral result,
but see \cite{gilletsoule2} for a rational result.

We define dual weight cohomology (with compact support) 
$H_W^i(X,\Q(n))^*$
to be the $i$th cohomology of the complex
$$CH^n(X_0)^*\leftarrow CH^n(X_1)^*\leftarrow CH^n(X_2)^*\leftarrow \cdots,$$
induced by contravariance of $CH^n$, and with 
$CH^n(X_i)^*$ placed in degree $2n+i$.
Note that this is the dual of the functor obtained via the 
contravariant analog of \cite[Thm.5.13]{jannsen} from the (contravariant) 
functor $CH^n(-)$ on the category $\C$.
We define dual weight cohomology of a field to be 
$$H^i_W(K,\Q(n))^*:= \colim_U H^i_W(U,\Q(n))^*,$$
where $U$ runs through smooth schemes with function field $K$.

\begin{lemma}\label{popoa}
We have $H^i_W(X,\Q(n))^*=0$ unless $2n\leq i\leq \dim X+n$. In particular,
$H^i_W(K,\Q(n))^*=0$ for every finitely generated field $K/k$ unless
$2n\leq i\leq \trdeg_k K+n$.
\end{lemma}

\proof This follows from the first property of
weight complexes together with $CH^n(T)=0$ for $n>\dim T$.
\proofend

It follows from Lemma \ref{popoa} that the niveau spectral sequence
\begin{equation}\label{uuiioo}
E^1_{s,t}=\bigoplus_{x\in X_{(s)}}H^{s+t}_W(k(x),\Q(n))^*\Rightarrow
H^{s+t}_W(X,\Q(n))^*
\end{equation}
is concentrated on and below the line $t=n$ and on and above the
line $s+t=2n$. If we let
$\tilde H^i_W(X,\Q(n))^*=E^2_{i-n,n}(X)$ be the homology of the complex
\begin{equation}\label{jannsenco}
\bigoplus_{x\in X_{(n)}}H^{2n}_W(k(x),\Q(n))^*\leftarrow
\cdots \leftarrow \bigoplus_{x\in X_{(d)}}H^{n+d}_W(k(x),\Q(n))^*,
\end{equation}
then we obtain a canonical and natural map
$$\gamma^*:\tilde  H^i_W(X,\Q(n))^*\to H^i_W(X,\Q(n))^* .$$

\subsection{Comparison}
We are going to check the hypothesis of \cite[Prop.5.16]{jannsen} 
to construct a functor between motivic cohomology with compact support
and weight cohomology. 
Recall that motivic cohomology with compact support is
defined as the cohomology of $C'(X)= \Hom_{D(Shv_{cdh})}(\Z^c(X),I^\cdot)$,
where $\rho^* \Z(n)\to I^\cdot$ is an injective resolution 
on the cdh-site. Then $C'$ is a covariant functor from the category of schemes
over $k$ with proper maps to the category of complexes with bounded above 
cohomology, which is contravariant for open embeddings.
Moreover, for proper $X$ we have $C'(X)=I^\cdot(X)$, and
a closed embedding $i:Y\to X$ with 
open complement $j:U\to X$ gives a short exact sequence
$$0\to C'(U)\to C'(X)\to C'(Y)\to 0.$$
Restricting $C'$ to smooth and proper $X$, we have 
$H^iC'(X)=0$ for $i>2n$, and a functorial isomorphism
$$H^{2n}C'(X)=H^{2n}I^\cdot(X)\cong \tau_{\geq 2n}I^\cdot(X) \cong CH^n(X).$$
by \eqref{smooth}. We obtain a morphism of functors on the category
of smooth and proper schemes,
$$C' = I^\cdot \to \tau_{\geq 2n}I^\cdot
\stackrel{\sim}{\longleftarrow}H^{2n}(I^\cdot)[-2n] = CH^n(-)[-2n]$$ 
Reversing all the arrows induced by arrows between schemes, but not 
by arrows between cohomology
theories in the proof of \cite[Prop.5.16]{jannsen} gives a natural 
transformation $H^i_c(X,\Z(n))\to H^i_W(X,\Z(n))$, 
hence a natural transformation
$$\pi^* :H^i_W(X,\Q(n))^*\to H^i_c(X,\Q(n))^* .$$

From now on we return to the situation $k$ finite.

\begin{proposition}\label{samefields}
Assume that every finitely generated field $K/k$ has a smooth and
projective model over $k$, and let $K$ be finitely generated 
of transcendence degree $d$ over $k$. 

a) The map $\pi^*$ induces isomorphisms
$$H^{n+d}_W(K,\Q(n))^*\stackrel{\sim}{\to} H^{n+d}_c(K,\Q(n))^*.$$
In particular, we have $\tilde H^i_W(X,\Q(n))^*\cong \tilde H^i_c(X,\Q(n))^*$.

b) If $d>n$, then $\pi^*$ induces isomorphisms
$$H^{n+d-1}_W(K,\Q(n))^*\stackrel{\sim}{\to} H^{n+d-1}_c(K,\Q(n))^*.$$ 
\end{proposition}

\proof
We proceed by induction on $d$.
Given $K$ of transcendence degree $d$, choose a smooth and projective
model $X$ of $K$ and compare \eqref{sseq} and \eqref{uuiioo}.

a) If $d<n$, then both terms vanish by Proposition \ref{vanbound}b) and 
Lemma \ref{popoa}. For $d=n$ we obtain 
$CH^n(X)\cong H^{n+d}_c(K,\Q(n))^*\cong H^{n+d}_W(K,\Q(n))^*$.
For $d>n$, we obtain
from $H^{n+d}_c(X,\Q(n))= H^{n+d}_W(X,\Q(n))=0$ a commutative diagram
with exact rows
$$\begin{CD}
\cdots @<<< \bigoplus_{x\in X_{(d-1)}}H^{n+d-1}_W(k(x),\Q(n))^*@<<<
H^{n+d}_W(K,\Q(n))^*@<<<  0\\
@EEE @| @VVV\\
\cdots @<<< \bigoplus_{x\in X_{(d-1)}}H^{n+d-1}_c(k(x),\Q(n))^*@<<<
H^{n+d}_c(K,\Q(n))^*@<<< 0.
\end{CD}$$

b) follows by a similar argument, noting that the $d_2$-differentials
originating from the terms in question end in terms considered in a), 
and there are no higher differentials.
\proofend

We obtain a commutative diagram
\begin{equation}\label{maind}
\begin{CD}
\tilde H^i_W(X,\Q(n))^* @= \tilde H^i_c(X,\Q(n))^*\\
@V\gamma^*VV @V\alpha^*VV  \\
H^i_W(X,\Q(n))^*@>\pi^*>>  H^i_c(X,\Q(n))^*.
\end{CD}
\end{equation}

\begin{proposition}
The following statements are equivalent:
\begin{enumerate}
\item Conjecture $P^n$.
\item The map $\pi^*$ is isomorphisms for all $X$.
\item We have $H^i_W(k,\Q(n))^*\cong H^i_c(k,\Q(n))^*$ for all $i$ and $k$.
\end{enumerate}
\end{proposition}

\proof
a) $\Leftrightarrow$ b): For smooth and proper $X$ this is clear.
In general, one does induction on the dimension and uses localization
sequences. 

b) $\Leftrightarrow$ c): One direction follows by taking colimits,
and the other by comparing the spectral sequences \eqref{sseq}
and \eqref{uuiioo}. 
\proofend

The following Proposition is analog to Proposition \ref{themapg}
and dual to \cite[Prop.3.4]{ichparshin}:

\begin{proposition}
The following statements are equivalent and follow from $\alpha^*$ being
an isomorphism:
\begin{enumerate}
\item For smooth and projective $X$, we have
$$\tilde H^i_W(X,\Q(n))^* \cong \begin{cases}
CH^n(X)^* &i=2n;\\
0 &\text{else}.
\end{cases}$$
\item The groups $H^i_W(k,\Q(n))^*$ vanish for $i\not=n+\trdeg k$.
\item The map $\gamma^*$ is an isomorphism for all $X$ and $i$.
\end{enumerate}
\end{proposition}

\proof
The proof is similar to Proposition \ref{themapg}. 

a) $\Rightarrow$ b): We proceed by induction on the transcendence
degree. Choose a smooth and projective model $X$ of $k$.
Since $H^i_W(X,\Q(n))$ is $CH^n(X)$ for $i=2n$ and
vanishes for $i\not=2n$, an inspection of
the spectral sequence \eqref{sseq} gives the result.

b) $\Rightarrow$ c) $\Rightarrow$ a) are clear. If $\alpha^*$ is 
an isomorphism, then so is $\pi^*$, and hence $\gamma^*$.

\proofend

\section{Beilinson's conjecture and duality}
Beilinson conjectured 
that over a finite field, rational and numerical
equivalence agrees up to torsion. This can be reformulated
to the following:

\medskip

\noindent{\bf Conjecture $D(n)$:}
{\it For all smooth and projective schemes $X$
over the finite field $\F_q$, the intersection pairing
gives a functorial isomorphism
$$ CH^n(X)_\Q\cong \Hom(CH_n(X),\Q).$$}

\medskip

Note that since both sides are countable, this implies finite
dimensionality. By the projection formula, the pairing induces 
a map of complexes

$$\begin{CD}
CH_n(X_0)@<<< CH_n(X_1)@<<< CH_n(X_2)@<<< \cdots\\
@VVV @VVV @VVV\\
CH^n(X_0)^*@<<< CH^n(X_1)^*@<<< CH^n(X_2)^* @<<< \cdots.
\end{CD}$$
Taking homology, we obtain a map
$$ \delta  :H_i^W(X,\Q(n)) \to H^i_W(X,\Q(n))^*.$$
Taking the limit over decreasing open sets with function field $K$,
$\delta$ induces a map $H_i^W(K,\Q(n)) \to H^i_W(K,\Q(n))^*$.
This in turn induces a map of complexes
$$\begin{CD}
\bigoplus_{x\in X_{(n)}}H_{2n}^W(k(x),\Q(n))
@<<< \bigoplus_{x\in X_{(n+1)}}H_{2n+1}^W(k(x),\Q(n))@<<< \cdots\\
@VVV @VVV \\
\bigoplus_{x\in X_{(n)}}H^{2n}_W(k(x),\Q(n))^*
@<<< \bigoplus_{x\in X_{(n+1)}}H^{2n+1}_W(k(x),\Q(n))^* @<<< \cdots,
\end{CD}$$
which gives the map $\tau$ making the following diagram commutative
\begin{equation}\label{alldiagram} \begin{CD}
H_i^c(X,\Q(n))@>\pi>>  H_i^W(X,\Q(n))@>\delta>>  H^i_W(X,\Q(n))^* 
@>\pi^* >>H^i_c(X,\Q(n))^*\\  \
@V\alpha VV @A\gamma AA @A\gamma^* AA @A\alpha^* AA \\
\tilde H_i^c(X,\Q(n))@>\beta >> \tilde H_i^W(X,\Q(n))@>\tau>> 
\tilde H^i_W(X,\Q(n))^* @= \tilde H^i_c(X,\Q(n))^*.
\end{CD}\end{equation}

\begin{lemma}
Conjecture $D(n)$ is equivalent to $\delta$ being an isomorphism for all 
$i$ and $X$, and implies that $\tau$ is an isomorphism for all $i$ and $X$.
\end{lemma}

\proof
The equivalence follows from the definition of $\delta$, and 
the statement about $\tau$ follows by a colimit argument.
\proofend

Parshin's conjecture and Beilinson's conjecture can be combined
into the following

\medskip

\noindent{\bf Conjecture $BP(n)$:}
{\it For all smooth and projective schemes $X$
over the finite field $\F_q$, the cup product pairing
$$ H^i(X,\Q(n))\times H^{2d-i}(X,\Q(d-n))\to \Q $$
is perfect.}

\medskip

\begin{proposition}\label{compareto}
For fixed $n$, the following statements are equivalent:
\begin{enumerate}
\item Conjecture $BP(n)$.
\item Conjectures $D(n)$, $P^n$ and $P_n$.
\item There are perfect pairings of finite
dimensional vector spaces
$$ H_i^c(X,\Q(n))\times  H^i_c(X,\Q(n))\to \Q $$
for all $X$, respectively smooth projective $X$. 
\item All maps in \eqref{alldiagram} are isomorphisms for all $X$,
respectively for all smooth and proper $X$.
\end{enumerate}
\end{proposition}

\proof
a) $\Leftrightarrow$ b): If $i>2n$, then the left hand side 
in $BP(n)$ vanishes, hence
perfectness is equivalent to the vanishing of 
$H^{2d-i}(X,\Q(d-n))\cong H_i^c(X,\Q(n))$
for $i>2n$, i.e. conjecture $P_n$ of \cite{ichparshin}. 
If $i<2n$, then the right hand side in $BP(n)$ vanishes, so perfectness
is equivalent to $P^n$. 
For $i=2n$, we recover conjecture $D(n)$. 

b) $\Leftrightarrow$ c):
Clearly conjecture $BP(n)$ is a special case of the assertion in c).
For the other direction, it suffices to construct a functorial map 
$H_i^c(X,\Q(n))\to  H^i_c(X,\Q(n))^*$ which is the intersection pairing
for smooth and projective $X$, and which is  
compatible with localization sequences on both sides. Indeed having
such a map one can use the usual devissage to reduce to the case that $X$
is smooth and projective. One way to construct such a map is 
to write $H_i^c(X,\Z(n))\cong \Hom_{DM^-}(\Z(n)[i],M^c(X))$,
$H^i_c(X,\Z(n))\cong \Hom_{DM^-}(M^c(X),\Z(n)[i])$, where 
$DM^-$ is Voevodsky's triangulated category of homotopy invariant
Nisnevich sheaves with transfers. Then the pairing is given by the
composition
\begin{multline*}
 \Hom_{DM^-}(\Z(n)[i],M^c(X))\times \Hom_{DM^-}(M^c(X),\Z(n)[i])
\to \\
\Hom_{DM^-}(\Z(n),\Z(n)) \cong \Hom_{DM^-}(\Z,\Z)\cong \Z,
\end{multline*}
using the cancellation theorem.

b) $\Leftrightarrow$ d)
Conjecture $P_n$, $D(n)$ and $P^n$ imply that the left square, middle 
horizontal maps, and right horizontal maps of \eqref{alldiagram} are 
isomorphisms for all $X$. 
Conversely, isomorphisms of the three upper maps of \eqref{alldiagram}
for smooth and proper $X$  
imply that $P_n$, $D(n)$, and $P^n$ hold, respectively. 
\proofend

\section{Parshin's conjecture and the four motivic theories}
Recall from \cite{friedvoe} that we have four motivic theories:
Motivic cohomology, motivic cohomology with compact support, 
motivic homology and motivic homology with compact support. 
All four theories are homotopy invariant and satisfy a 
projective bundle formula.
Motivic cohomology is contravariant, has a Mayer-Vietoris long
exact sequence for Zarsiki covers, and a long exact sequence
for abstract blow-ups. Motivic cohomology is contravariant
for proper maps, covariant for quasi-finite flat maps, and
satisfies a localization long exact sequence (which implies
in particular Mayer-Vietoris and abstract blow-up 
long exact sequences). Motivic homology and motivic homology
with compact support satisfy the dual properties.
The theories are related by the following diagram
$$ \begin{CD}
H^i_c(X,\Q(n))@>proper>\sim>  H^i(X,\Q(n))\\
@VsmoothV\cong V @VsmoothV\cong V \\
H_j(X,\Q(m))@>proper>\sim> H_j^c(X,\Q(m))
\end{CD}$$
The horizontal maps are isomorphisms for proper $X$, and the 
vertical maps are isomorphisms if $X$ is smooth of pure dimension $d$,
and $m+n=d$ and $j+i=2d$. The functorialities suggest
that groups diagonally opposite should be in some form of duality; 
we saw that with rational coefficients, this is equivalent to deep 
conjectures, for a result with torsion coefficients see \cite{ichdual}.

The following diagram describes the range where these groups
can be non-zero, where they can be non-zero 
assuming Parshin's conjecture, where they can be non-zero
assuming Parshin's conjecture
plus smoothness of $X$, and where they can be non-zero 
assuming Parshin's conjecture plus properness of $X$, respectively.
The bold faced inequalities indicate that they are strong enough
to recover Parshin's conjecture.

\medskip

\noindent
{\renewcommand\arraystretch{1.5}
\begin{tabular}{l|c|c|c|c}
&Coh compact sup&Mot Cohomology&Mot Homology&Borel-Moore hom\\
\hline
&$H^i_c(X,\Q(n))$&$ H^i(X,\Q(n))$ &$ H_j(X,\Q(m))$ &$ H_j^c(X,\Q(m))$ \\
\hline
always&$i\leq n+d$ &$i\leq n+d $&$j\geq m$&$j\geq 2m $\\
& &$i\leq 2n$ X smooth & $j\geq 2m$ X proper \\
Parshin $\Rightarrow$ &
${\bf 2n\leq i}\leq n+d$ &$ n\leq i\leq n+d $&$ m\leq j\leq m+d$
&$2m\leq j\leq m+d$ \\
P+smooth&  &$n\leq i\leq 2n$&$m\leq {\bf j\leq 2m}$&\\
P+proper&  &${\bf 2n\leq i}\leq n+d$&$2m\leq j\leq m+d$&
\end{tabular}}

\bigskip 

\proof
The first row follows from the definitions (and that the cdh-cohomological
dimension agress with the dimension).
Since Borel-Moore homology $H_j^c(X,\Q(m))$ is isomorphic to
higher Chow groups $CH_m(X,j-2m)$, they can only be non-zero for 
$j\geq 2m$. The second row is the translation of this fact into a statement
for motivic cohomology for smooth $X$, and for motivic homology
for proper $X$.

The results under Parshin's conjecture for Borel-Moore homology
and motivic cohomology with compact support can be obtained by using 
induction on the dimension and the localization sequences. 
To obtain them for motivic homology and cohomology, 
one uses the isomorphisms 
$H_i(X,\Z(n))\cong H^{2d-i}_c(X,\Z(d-n))$
and $H^i(X,\Z(n))\cong H_{2d-i}^c(X,\Z(d-n))$
for a smooth scheme $X$ of dimension $d$ to obtain the result for 
smooth schemes. Then induction on the dimension and the blow-up long exact 
sequences gives results for all schemes.

The extra information for the smooth and proper case in case
of homology and cohomology is obtained by comparing to the other
theories.
\proofend


The bold faced inequalities were a
motivation to write this paper: It might be difficult to prove
a statement which only holds for smooth and proper $X$, as in 
the case of higher Chow groups. It might be easier to prove
a statement which holds for all smooth schemes (motivic homology), or 
all proper schemes (motivic cohomology), or all schemes
(motivic cohomology with compact support).

\end{document}